\documentclass{llncs}

\usepackage{amsmath}
\usepackage{amssymb}
\usepackage{graphicx}
\usepackage[colorinlistoftodos]{todonotes}
\usepackage[colorlinks=true, allcolors=blue]{hyperref}
\usepackage{setspace}
\usepackage[margin=1.5in]{geometry}
\usepackage{lastpage}
\usepackage{fancyhdr}
\begin{document}
\cfoot{\thepage}
\title{On the overestimation of the largest eigenvalue of a covariance matrix}
\titlerunning{Hamiltonian Mechanics}  % abbreviated title (for running head)
\author{Soufiane Hayou}
\institute{{\it Department of Applied Mathematics, Ecole polytechnique (Paris)}\\
\email{soufiane.hayou@polytechnique.edu}}

\maketitle              % typeset the title of the contribution

\begin{abstract}
In this paper, we use a new approach to prove that the largest eigenvalue of the sample covariance matrix of a normally distributed vector is bigger than the true largest eigenvalue with probability 1 when the dimension is infinite.
\end{abstract}

\section{Introduction}

It is well known that when the number of samples is big compared to the dimension of the variables, one can accurately estimate the correlation matrix using the sample covariance matrix :\\
Let $p$ be the dimension, $n$ the number of samples and $(X_i)_{1\leq i \leq T}$ the observations. The sample covariance matrix is defined by :
\begin{equation*}
S = \frac{1}{n} \sum_{i=1}^n (X_i - \bar{X}) (X_i - \bar{X})^t
\end{equation*}
where $X^t$ is the transpose of $X$, and $\bar{X}$ is the empirical mean.\\

However, the sample covariance matrix is not a good estimator when $p$ is of the same range as $n$. Actually, we empirically observe that the sample covariance matrix tends to overestimate (underestimate) the largest (smallest) eigenvalues for a large class of covariance matrices.\\ 

There is an extensive amount of research papers on this topic. We find (in general) three main categories of the proposed approaches : Shrinkage to a target (e.g. Ledoit and Wolf \cite{ledoit}), Random Matrix Theory (N. El Karoui \cite{elkaroui}, Bouchaud and Bun \cite{rmt}) and Optimization under constraints (e.g. under a constraint on the condition number like in \cite{condition}). In particular, it has been proved in \cite{karoui2}, under some conditions on the limiting spectrum (limiting density), that the largest sample eigenvalue has a Tracy-Widom distribution when $p$ goes to infinity (with $\frac{p}{n}$ uniformly bounded). Using this, one can prove easily that in this setting, with probability 1, we overestimate the largest eigenvalue. However, these conditions are sometimes hard to verify, and the sufficient conditions in \cite{karoui2} impose on the empirical density to have a limiting density. Our goal is to overcome this problem and give some easily-verified conditions on the spectrum without any restrictions on the limit.\\

In this paper, we use a result from \cite{muirhead} on the upper bound of the probability distribution of the largest eigenvalue to show that, under some assumptions, the probability of the event $\{$largest eigenvalue of the sample $>$ largest eigenvalue of the true covariance matrix$\}$ converges to 1 when the dimension goes to infinity.\\

In what follows, S is the sample covariance matrix, $\Sigma$ is the true covariance matrix, $l_1 \geq l_2 \geq ... \geq l_p$ the eigenvalues of S, $\lambda_1 \geq \lambda_2 \geq ... \geq \lambda_p$ the eigenvalues of $\Sigma$ and $q = \frac{p}{n} < 1$ is the ratio of the dimension over the sample size. In the first section, we present some classical results on Wishart matrices (a generalization to multiple dimensions of the chi-squared distribution) from Random Matrix Theory, we show the main results in sections 2 and 3.\\

\section{Wishart matrices and the distribution of the eigenvalues}
In this section, we present the Wishart distribution and some related results from random matrix theory. We start by recalling the definition of a Wishart matrix.\\

{\bf Definition} : A $p \times p$ matrix M is said to have a Wishart distribution with covariance matrix $\Sigma$ and degrees of freedom n if there exists $X \sim N_{n \times p}(\mu, \Sigma)$ such that $M = X^t X$. We denote it by $M \sim W_p(n, \Sigma)$.\\

When $n \geq p$, the Wishart distribution has a density function given by : 
\begin{equation}
f(M) = \frac{2^{-np/2}}{\Gamma_p(n/2) (det(\Sigma))^{n/2}} etr(-\frac{1}{2} \Sigma^{-1} M) (det M)^{(n-p-1)/2}
\end{equation}
where $etr$ is the exponential of the trace, $\Gamma_p$ is the generalized gamma function.\\

When $X \sim N_{n \times p}(\mu, \Sigma)$, the sample covariance matrix $S = \frac{1}{n} X X^t$ has the Wishart distribution $W_p(n-1, \frac{1}{n}\Sigma)$ (see \cite{bejan} for the proof).\\

\subsection{Joint distribution of the eigenvalues}
Let $M \sim W_p(n,\Sigma)$ where $n >p$, then the joint distribution of the eigenvalues $l_1 \geq l_2 \geq ... \geq l_p$ is given by (see \cite{bejan}): 
\begin{equation}
g(l_1, l_2, ..., l_p) = \\
\frac{\pi^{p^2/2} \times 2^{-np/2} (det \Sigma)^{-n/2}}{\Gamma_p(n/2) \Gamma_p(p/2)} \prod_{i=1}^p l_i^{(n-p-1)/2} \prod_{j>i}^p (l_i - l_j) \int_{O_p} etr(-\frac{1}{2} \Sigma^{-1} HLH^t) (dH)
\end{equation}

where $L = diag(l_1, l_2, ..., l_p)$, and the integral is over the orthogonal group $O_p$ with respect to the Haar measure (see \cite{Elizabeth}).\\
In general, the integral is hard to estimate. However, when $\Sigma = \lambda I$, we have : 
\begin{align*}
\int_{O_p} etr(-\frac{1}{2} \Sigma^{-1} HLH^t) (dH) &= \int_{O_p} etr(-\frac{1}{2 \lambda} HLH^t) (dH)\\
&= etr(-\frac{1}{2 \lambda}L) \int_{O_p} (dH)\\
&= \exp(-\frac{1}{2 \lambda} \sum_{i=1}^p l_i)
\end{align*}

The Haar measure is invariant by rotation, that means for any orthogonal matrix $Q$, one has :
\begin{equation*}
d(QH) = dH
\end{equation*}

Using this and the fact that there exists an orthogonal matrix $Q$ such that $\Sigma^{-1} = Q D^{-1} Q^t$ where $D = diag(\lambda_1, \lambda_2, ..., \lambda_p)$, we can prove that the previous distribution depends only on the eigenvalues of $\Sigma$.\\

We know that $S \sim W_p(n-1, \frac{1}{n}\Sigma)$, so the joint distribution of the eigenvalues $l_1 \geq l_2 \geq ... \geq l_p$ of the sample covariance matrix is given by : 
\begin{equation*}
g(l_1, l_2, ..., l_p) = \\
\frac{\pi^{p^2/2} (det \Sigma)^{-(n-1)/2} }{ \Gamma_p((n-1)/2) \Gamma_p(p/2)} (\frac{n-1}{2})^{\frac{p (n-1)}{2}} \prod_{i=1}^p l_i^{(n-p-2)/2} \prod_{j>i}^p (l_i - l_j) \int_{O_p} etr(-\frac{1}{2} n\Sigma^{-1} HLH^t) (dH)
\end{equation*}

\subsection{Distribution of the largest eigenvalue of the sample covariance matrix}
The cumulative distribution function of the largest eigenvalue of S is given by :
\begin{equation}
\mathbb{P}(l_1 < x) = \frac{\Gamma_p(\frac{p+1}{2})}{\Gamma_p(\frac{p+n}{2})} det(\frac{n-1}{2} \Sigma^{-1})^{(n-1)/2} F_{1,1}(\frac{n-1}{2}; \frac{n+p}{2}; -\frac{n}{2} x \Sigma^{-1})
\end{equation}
where $F_{1,1}$ is hypergeometric function with a matrix argument (see \cite{eduardo}). This function is hard to evaluate, which makes the previous formula hard to use directly. The next result was first proved by R.J. Murhead in \cite{muirhead}.\\

{\bf Theorem I (Muirhead) }: Let x be a nonnegative real number. The following inequalities hold for any $p$ and $n$ such that $p<n$ :
\begin{equation}
\mathbb{P}(l_1 \leq x) \leq \prod_{i=1}^p \mathbb{P}(\chi_n^2 \leq \frac{nx}{\lambda_i})
\end{equation}
\begin{equation}
\mathbb{P}(l_p \leq x) \geq 1 - \prod_{i=1}^p \mathbb{P}(\chi_n^2 \geq \frac{nx}{\lambda_i})
\end{equation}

where $\chi_n^2$ is a chi-square random variable with $n$ degrees of freedom. \\

Since finding bounds on hypergeometric functions is still a subject of interest, we cannot fully check the quality of this bounds (The only way to check it, is by simulation).

\subsection{Special case : Marchenko-Pastur distribution}
When $\Sigma = I$, the Marchenko-Pastur theorem states that the empirical distribution of the sample eigenvalues converges when $p \rightarrow \infty$ (with $q=\frac{p}{n}$ fixed) to the Marchenko-Pastur distribution given by :
\begin{equation*}
mp(x) = \frac{1}{2 \pi} \frac{\sqrt{(\lambda_+ - x)(x - \lambda_-)}}{qx} \hspace*{0.1cm} 1_{\lambda_- \leq x \leq \lambda_+}
\end{equation*}
where $\lambda_+ = (1 + \sqrt{q})^2$ and $\lambda_- = (1 - \sqrt{q})^2$.\\

Figure 1 shows the Marchenko-Pastur distribution for q=0.1 .

\begin{figure*}
  \centerline{\includegraphics[width=0.6\linewidth]{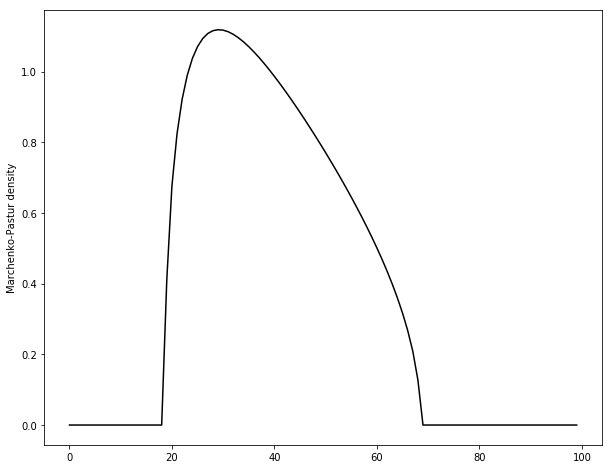}}
  \caption{Marchenko-Pastur distribution for q = 0.1}
\end{figure*}

\section{Overestimating the largest eigenvalue}
In this section, we prove that the probability of the event {$l_1 > \lambda_1$} converges to 1 when the dimension goes to infinity (with $q = \frac{p}{n}$ fixed) under some constraints.\\
The following lemmas will be proved in order to be used in the proof of the result.\\

{\bf Lemma 1 }: Let $(\Omega, \mathbb{P}, \mathbb{F})$ be a probability space, and $(A_n)$, $(B_n)$ two sequences of events (not necessary independent), with $\lim_{n \rightarrow \infty} \mathbb{P}(B_n) = 1$, then we have : \\
\begin{equation}
\limsup_{n \rightarrow \infty} \mathbb{P}(A_n \cap B_n) = \limsup_{n \rightarrow \infty} \mathbb{P}(A_n)
\end{equation}

{\it Proof :}\\
We have,

\begin{equation*}
\mathbb{P}(A_n \cap B_n) = \mathbb{P}(A_n) + \mathbb{P}( B_n) - \mathbb{P}(A_n \cup B_n)
\end{equation*}
and using the fact that $\mathbb{P}(A_n \cup B_n) \geq \mathbb{P}(B_n)$ and $\mathbb{P}(B_n) \rightarrow 1$, we conclude.\\

{\bf Lemma 2 }: Let $\chi_n^2$ be a chi-square random variable with $n$ degrees of freedom, $(a_n)$ is a serie of positive real numbers. Then for any increasing and continuous function $f$ and $\epsilon > 0$, we have :
\begin{equation}
\limsup_{n \rightarrow \infty} f(\mathbb{P}(N(0,1) \leq \sqrt{n} (a_n -1) - \epsilon)) \leq \limsup_{n \rightarrow \infty} f(\mathbb{P}(\chi_n^2 \leq n a_n)) \leq \limsup_{n \rightarrow \infty} f(\mathbb{P}(N(0,1) \leq \sqrt{n} (a_n -1) + \epsilon))
\end{equation}\\

and,\\
\begin{equation}
\limsup_{n \rightarrow \infty} f(\mathbb{P}(N(0,1) \geq \sqrt{n} (a_n -1) + \epsilon)) \leq \limsup_{n \rightarrow \infty} f(\mathbb{P}(\chi_n^2 \geq n a_n)) \leq \limsup_{n \rightarrow \infty} f(\mathbb{P}(N(0,1) \geq \sqrt{n} (a_n -1) - \epsilon))
\end{equation}\\

{\it Proof :}\\

We know that $\chi_n^2 =^{d} Z_1^2 + Z_2^2 + ... + Z_n^2$ where $(Z_i)_{1 \leq i \leq n}$ are standard normal variables (the equality holds in distribution). Since Gaussian variables have moments of any order, we can use the Central Limit Theorem and we have : 
\begin{equation*}
\sqrt{n} (\frac{\chi_n^2}{n} - 1) \rightarrow_{d} N(0,1)
\end{equation*}
We can write this as :
\begin{equation}
\frac{\chi_n^2}{n} =_{d} 1 + \frac{1}{\sqrt{n}} Z + \frac{1}{\sqrt{n}} \epsilon_n
\end{equation}
where $Z \sim N(0,1)$ and $\epsilon_n = o_d(1)$ ($o_d$ means $\epsilon_n$ converges to 0 in distribution). We know that the convergence in distribution to a constant implies the convergence in probability to the same constant. That means, we can write $\epsilon_n = o_{\mathbb{P}}(1)$.\\

Now let $\epsilon>0$. We have : 

\begin{align*}
\mathbb{P}(\chi_n^2 \leq n a_n) &= \mathbb{P}(Z + \epsilon_n \leq \sqrt{n} (a_n-1))\\
&= \mathbb{P}(Z + \epsilon_n \leq \sqrt{n} (a_n-1) \hspace{0.1cm} , \hspace{0.1cm} |\epsilon_n|<\epsilon) \\
&+ \mathbb{P}(Z + \epsilon_n \leq \sqrt{n} (a_n-1) \hspace{0.1cm} | \hspace{0.1cm} |\epsilon_n|>\epsilon) \mathbb{P}(|\epsilon_n|>\epsilon)\\
\end{align*}
We also have that,
\begin{center}
$\mathbb{P}(Z + \epsilon \leq \sqrt{n} (a_n-1), \hspace{0.1cm} |\epsilon_n|<\epsilon) \leq \mathbb{P}(Z + \epsilon_n \leq \sqrt{n} (a_n-1) \hspace{0.1cm} , \hspace{0.1cm} |\epsilon_n|<\epsilon)$\\
$\mathbb{P}(Z + \epsilon_n \leq \sqrt{n} (a_n-1) \hspace{0.1cm} , \hspace{0.1cm} |\epsilon_n|<\epsilon) \leq \mathbb{P}(Z -\epsilon \leq \sqrt{n} (a_n-1), \hspace{0.1cm} |\epsilon_n|<\epsilon)$ 
\end{center}

Now we use Lemma 1 with the sequences $A_n = \{Z \leq \sqrt{n} (a_n-1) - \epsilon\}$ and $B_n = \{ |\epsilon_n| < \epsilon\} $ for the left-hand inequality and the sequences $A_n = \{Z \leq \sqrt{n} (a_n-1) + \epsilon\}$ and $B_n = \{ |\epsilon_n| < \epsilon\} $ for the right-hand inequality. We conclude since $\lim_{n \rightarrow \infty} \mathbb{P}(|\epsilon_n|>\epsilon) = 0$ and f is increasing and continuous. \\

The second inequality can be deduced from the first one using the function $x \rightarrow g(x) = - f(1 - x)$ which is increasing and continuous. \\

{\bf Lemma 3 }: Let F be the cumulative distribution function of a standard normal variable. Then the following inequality holds for any real number x : 
\begin{equation}
\frac{1}{x + \sqrt{x^2+4}} \leq \sqrt{\frac{\pi}{2}} e^{\frac{x^2}{2}} (1 - F(x)) \leq \frac{1}{x + \sqrt{x^2+\frac{8}{\pi}}}
\end{equation}

{\it Proof :}\\

We use the following inequality (Formula 7.1.13 from \cite{abramowitz}), for any real number y :
\begin{equation}
\frac{1}{y + \sqrt{y^2 + 2}} \leq e^{y^2} \int_y^{\infty} e^{-t^2} dt \leq \frac{1}{y + \sqrt{y^2 + \frac{4}{\pi}}}
\end{equation}
We use the new variable $x$ defined by $x = \sqrt{2} y$, then $\int_y^{\infty} e^{-t^2} dt = \frac{1}{\sqrt{2}} \int_x^{\infty} e^{-\frac{t^2}{2}} dt $. The inequality follows from that.\\

Now, we prove the main result of this paper. \\

{\bf Theorem II (main result) }: Let $p,n$ be two positive integers such that $q = \frac{p}{n}<1$ is fixed, $(\lambda_{1,p} \geq \lambda_{2,p} \geq ... \geq \lambda_{p,p})_{p>0}$ a sequence of spectrums (spectrum of $\Sigma$) and $(l_{1,p} \geq l_{2,p} \geq ... \geq l_{p,p})_{p>0}$ the sequence of the corresponding sample spectrums (spectrum of $S$). 
For any p, we define the set $J_p$ by :
\begin{equation}
J_p = \{i : for \hspace{0.1cm} all \hspace{0.1cm} m \geq p, \hspace{0.1cm}  |\frac{\lambda_{1,m}}{\lambda_{i,m}} - 1| < \frac{1}{\sqrt{m}} \}
\end{equation}
and the cardinal number of $J_p$ (number of elements in $J_p$):
\begin{equation*}
\phi(p) = |J_p|
\end{equation*}
Then, $\phi(p)$ is increasing and has a limit when $p \rightarrow \infty$, and there exists a constant $c>0$ such that for any $p>0$:
\begin{equation}
\limsup_{p \rightarrow \infty} \mathbb{P}(l_{1,p} \leq \lambda_{1,p}) \leq \limsup_{p \rightarrow \infty} e^{- c \times \phi(p)} = e^{- c \times \lim_{p \rightarrow \infty} \phi(p)}
\end{equation}
As a result, we have :
\begin{equation}
\lim_{p \rightarrow \infty} \phi(p) = \infty  \Rightarrow \lim_{p \rightarrow \infty} \mathbb{P}(l_{1,p} \leq \lambda_{1,p}) = 0
\end{equation}
The previous result can be interpreted as follows : When the largest eigenvalue is not isolated from the spectrum in the sense that the set $J_p$ has infinite number of elements when p goes to infinity, then, with probability 1, the sample covariance matrix overestimates it when the dimension is infinite (with $q$ fixed). \\

{\it Proof :}\\

Let x be a nonnegative real number and $\epsilon>0$. We recall Muirhead's upper bound : 
\begin{equation*}
\mathbb{P}(l_{1,p} \leq x) \leq \prod_{i=1}^p \mathbb{P}(\chi_n^2 \leq \frac{nx}{\lambda_{i,p}})
\end{equation*}

We want to prove that the right-hand quantity converges to 0 when $p \rightarrow \infty$. Since it is nonnegative, we will prove that the limit superior of the same quantity converges to 0. We use the following notations :\\

\begin{center}
$x_{i,n} = \sqrt{n}(\frac{\lambda_{1,p}}{\lambda_{i,p}} -1)$, $i \leq p = qn$ \\
\end{center}

\[ a_{i,n} = \begin{cases} 
      \mathbb{P}(\chi_n^2 \leq \frac{n \lambda_{1,p}}{\lambda_{i,p}}) & i \leq p=q n \\
      1 & i > p = q n 
   \end{cases}
\]\\

and 

\[ b_{i,n} = \begin{cases} 
      \mathbb{P}(N(0,1) \leq \sqrt{n} (\frac{\lambda_{1,p}}{\lambda_{i,p}} -1) + \epsilon) & i \leq p=q n \\
      1 & i > p = q n 
   \end{cases}
\]\\

Then we have :
\begin{align*}
\limsup_{p \rightarrow \infty} \log(\prod_{i=1}^p \mathbb{P}(\chi_{p/q}^2 \leq \frac{p/q \lambda_{1,p}}{\lambda_{i,p}})) & \leq \limsup_{n \rightarrow \infty} \sum_{i=1}^{\infty} \log(a_{i,n})\\
& \leq  \sum_{i=1}^{\infty} \limsup_{n \rightarrow \infty} \log(a_{i,n})\\
& \leq  \sum_{i=1}^{\infty} \limsup_{n \rightarrow \infty} \log(b_{i,n}) \hspace{1cm} \textnormal{(Lemma 2)}\\
& \leq  \sum_{i=1}^{\infty} \limsup_{n \rightarrow \infty} \log(1 - z_{i,n}) \hspace{1cm} \textnormal{(Lemma 3)}\\
\end{align*}
where, for $i \leq qn$,
\begin{equation}
z_{i,n} = \sqrt{\frac{2}{\pi}} e^{-\frac{(x_{i,n} + \epsilon)^2}{2}} \frac{1}{x_{i,n}+\epsilon + \sqrt{(x_{i,n} + \epsilon)^2 + 4}}
\end{equation}
and $z_{i,n} = 0$ otherwise.
We know that for any real number $x<1$, we have $\log(1 - x) \leq - x$, and it is clear that $z_{i,n}<1$ for any i and n (since $x_{i,n}, \epsilon \geq 0$), therefore, for any p (and $n=\frac{p}{q}$) we have : 

\begin{align*}
\sum_{i=1}^{\infty} \sup_{m \geq n} \log(1 - z_{i,m})
& \leq \sum_{i=1}^{\infty} \sup_{m \geq n} - z_{i,m}\\
& \leq \sum_{i=1}^{\infty} - \inf_{m \geq n} z_{i,m}\\
& = - \sum_{i=1}^{\infty} \inf_{m \geq n} z_{i,m}\\
\end{align*}

Now, using the inequality $(x_{i,m}+\epsilon)^2 \leq 2(x_{i,m}^2+\epsilon^2)$, and the fact that for i $\in$ $J_p = J_{qn}$, $|x_{i,m}| \leq 1$ ($ \forall m \geq n$), we have that $\forall i \in J_p$ , $\forall m \geq n$ : 
\begin{align*}
z_{i,m} & \geq \sqrt{\frac{2}{\pi}} e^{-\frac{2 (x_{i,m}^2 + \epsilon^2)}{2}} \frac{1}{x_{i,m}+\epsilon + \sqrt{(x_{i,m} + \epsilon)^2 + 4}}\\
\textnormal{so that,} \\
z_{i,m} & \geq \sqrt{\frac{2}{\pi}} e^{-\frac{2 (1 + \epsilon^2)}{2}} \frac{1}{1+\epsilon + \sqrt{(1 + \epsilon)^2 + 4}}\\
\textnormal{therefore,} \\
\inf_{m \geq n} z_{i,m} & \geq \sqrt{\frac{2}{\pi}} e^{-\frac{2 (1 + \epsilon^2)}{2}} \frac{1}{1+\epsilon + \sqrt{(1 + \epsilon)^2 + 4}}\\
\end{align*}

That gives us the following inequality : 
\begin{equation*}
\sum_{i=1}^{\infty} \inf_{m \geq n} z_{i,m} \geq \phi(qn) \sqrt{\frac{2}{\pi}} e^{-\frac{2 (1 + \epsilon^2)}{2}} \frac{1}{1+\epsilon + \sqrt{(1 + \epsilon)^2 + 4}}
\end{equation*}

Note that $\phi(p)$ is increasing, so it has a limit. Using monotone convergence theorem (since $\forall i,m \log(1 - z_{i,m}) \leq 0$),
\begin{align*}
\sum_{i=1}^{\infty} \limsup_{n \rightarrow \infty} \log(1 - z_{i,n}) &= \lim_{n \rightarrow \infty} \sum_{i=1}^{\infty} \sup_{m \geq n} \log(1 - z_{i,m}) \\
& \leq - \lim_{n \rightarrow \infty} \sum_{i=1}^{\infty} \inf_{m \geq n} z_{i,m}\\
& \leq - \lim_{n \rightarrow \infty} \phi(qn) c_{\epsilon}
\end{align*}

Therefore,
\begin{equation*}
\limsup_{p \rightarrow \infty} \prod_{i=1}^p \mathbb{P}(\chi_{p/q}^2 \leq \frac{p/q \lambda_{1,p}}{\lambda_{i,p}}) \leq  e^{- c_{\epsilon} \times \lim_{p \rightarrow \infty} \phi(p)}
\end{equation*}
where $c_{\epsilon} = \sqrt{\frac{2}{\pi}} e^{-\frac{2 (1 + \epsilon^2)}{2}} \frac{1}{1+\epsilon + \sqrt{(1 + \epsilon)^2 + 4}}$.\\

Since this is true for any $\epsilon>0$, then (in both cases $\lim_{p \rightarrow \infty} \phi(p)$ finite or infinite) we have that : 
\begin{equation}
\limsup_{p \rightarrow \infty} \prod_{i=1}^p \mathbb{P}(\chi_{p/q}^2 \leq \frac{p/q \lambda_{1,p}}{\lambda_{i,p}}) \leq  e^{- c \times \lim_{p \rightarrow \infty} \phi(p)}
\end{equation}
where $c = \sqrt{\frac{2}{\pi}} e^{-1} \frac{1}{1 + \sqrt{5}}$.\\

Using Muirhead's inequality we conclude that : 
\begin{equation}
\limsup_{p \rightarrow \infty} \mathbb{P}(l_{1,p} \leq \lambda_{1,p}) \leq e^{- c \times \lim_{p \rightarrow \infty} \phi(p)}
\end{equation}\\

Note that this theorem is valid for any sequence $(x_p)$ (and not just the the sequence ($\lambda_{1,p}$)). This proves the following corollary :\\

{\bf Corollary 1 }: Let $(x_p)_{p>0}$ be a sequence of positive real numbers. We define the set $J_p(x)$ by :
\begin{equation}
J_p(x) = \{i : for \hspace{0.1cm} all \hspace{0.1cm} m \geq p, \hspace{0.1cm}  |\frac{x_m}{\lambda_{i,m}} - 1| < \frac{1}{\sqrt{m}} \}
\end{equation}
and,
\begin{equation}
\phi(x,p) = |J_p(x)|
\end{equation}

Then for any sequence x such that $\lim_{p \rightarrow \infty} \phi(x,p) = \infty$ we have that $\lim_{p \rightarrow \infty} \mathbb{P}(l_{1,p} \leq x_p) = 0$. \\
This means that for any sequence x that is not isolated from the spectrum when the dimension goes to infinity, with probability 1, the sequence of the largest sample eigenvalues is greater than $x$ (element-wise). Note that if $x$ has a limit $\gamma$, then the largest sample eigenvalue is bigger than $\gamma$ with probability 1.

\section{Underestimating the smallest eigenvalue}
In this section, we show that the sample covariance matrix underestimates the smallest eigenvalue with a probability 1 when the dimension is infinite. The constraints for this result to be true are slightly different from the previous theorem, but the proof is similar.\\

{\bf Theorem III }: Let $p,n$ be two positive integers such that $q = \frac{p}{n}<1$ is fixed, $(\lambda_{1,p} \geq \lambda_{2,p} \geq ... \geq \lambda_{p,p})_{p>0}$ a sequence of spectrums (spectrum of $\Sigma$) and $(l_{1,p} \geq l_{2,p} \geq ... \geq l_{p,p})_{p>0}$ the sequence of the corresponding sample spectrums (spectrum of $S$).  
For any positive integer p and positive real number $\kappa$, we define the set $H_{p,\kappa}$ by :
\begin{equation}
H_{p, \kappa} = \{i : for \hspace{0.1cm} all \hspace{0.1cm} m \geq p, \hspace{0.1cm}  |\frac{\lambda_{m,m}}{\lambda_{i,m}} - 1| < \frac{\kappa}{\sqrt{m}} \}
\end{equation}
and the cardinal number of the set $H_{p, \kappa}$:
\begin{equation*}
\xi(\kappa, p) = |H_{p, \kappa}|
\end{equation*}
Then, for any $\kappa < \sqrt{\frac{2}{\pi}}$ there exists $c_{\kappa} > 0$ such that :
\begin{equation}
\liminf_{p \rightarrow \infty} \mathbb{P}(l_{p,p} \leq \lambda_{p,p}) \geq  1 - \limsup_{p \rightarrow \infty} e^{- c_{\kappa} \xi(\kappa,p)}
\end{equation}
As a result, we have :\\
for any $\kappa < \sqrt{\frac{2}{\pi}}$,
\begin{equation}
\lim_{p \rightarrow \infty} \xi(\kappa, p) = \infty \Rightarrow \lim_{p \rightarrow \infty}\mathbb{P}(l_{p,p} \leq \lambda_{p,p}) = 1
\end{equation} \\

{\it Proof }:\\

Let $\epsilon>0$ and $\kappa>0$ such  that $\kappa < \sqrt{\frac{2}{\pi}}$.
We recall Muirhead's lower bound for the distribution of the smallest eigenvalue :
\begin{equation*}
\mathbb{P}(l_{p,p} \leq x) \geq 1 - \prod_{i=1}^p \mathbb{P}(\chi_n^2 \geq \frac{nx}{\lambda_{i,p}})
\end{equation*}

Similar to the previous proof, we use the following notations :\\

\begin{center}
$x_{i,n} = \sqrt{n}(\frac{\lambda_{p,p}}{\lambda_{i,p}} -1), \hspace*{0.5cm} i\leq p=qn$ \\
\end{center}

\[ a_{i,n} = \begin{cases} 
      \mathbb{P}(\chi_n^2 \geq \frac{n \lambda_{p,p}}{\lambda_{i,p}}) & i \leq p=q n \\
      1 & i > p = q n 
   \end{cases}
\]\\

and 

\[ b_{i,n} = \begin{cases} 
      \mathbb{P}(N(0,1) \geq \sqrt{n} (\frac{\lambda_{p,p}}{\lambda_{i,p}} -1) - \epsilon) & i \leq p=q n \\
      1 & i > p = q n 
   \end{cases}
\]\\

Then we have that,\\
\begin{align*}
\limsup_{p \rightarrow \infty} \log(\prod_{i=1}^p \mathbb{P}(\chi_{p/q}^2 \geq \frac{p/q \lambda_{1,p}}{\lambda_{i,p}})) & \leq \limsup_{n \rightarrow \infty} \sum_{i=1}^{\infty} \log(a_{i,n})\\
& \leq  \sum_{i=1}^{\infty} \limsup_{n \rightarrow \infty} \log(a_{i,n})\\
& \leq  \sum_{i=1}^{\infty} \limsup_{n \rightarrow \infty} \log(b_{i,n}) \hspace{1cm} \textnormal{(Lemma 2)}\\
\end{align*}

Since $\forall i,m, \hspace*{0.1cm} 0 < b_{i,m} \leq 1$, we have that :
\begin{align*}
\sum_{i=1}^{\infty} \sup_{m \geq n} \log(b_{i,m}) & \leq \sum_{i \in H_{qn,\kappa}} \sup_{m \geq n} \log(b_{i,m}) \\
& \leq \sum_{i \in H_{qn,\kappa}} \sup_{m \geq n} \log(z_{i,m})
\end{align*}

where for $i \leq qm$,
\begin{equation*}
z_{i,m} = \sqrt{\frac{2}{\pi}} e^{-\frac{(x_{i,m} - \epsilon)^2}{2}} \frac{1}{x_{i,m}-\epsilon + \sqrt{(x_{i,m} - \epsilon)^2 + \frac{8}{\pi}}}
\end{equation*}
and $z_{i,m} = 1$ otherwise.\\

Now let n>0, we have that $\forall i \in H_{qn,\kappa}, m\geq n$ : 
\begin{equation*}
\frac{1}{x_{i,m} - \epsilon + \sqrt{(x_{i,m} - \epsilon)^2 + \frac{8}{\pi}}} \leq \frac{1}{-\kappa - \epsilon + \sqrt{\epsilon^2 + \frac{8}{\pi}}}  
\end{equation*}
so that,
\begin{equation*}
\sup_{m \geq n} \log(z_{i,m}) \leq c_{\epsilon} 
\end{equation*}
where $c_{\epsilon} = \log(\sqrt{\frac{2}{\pi}} \frac{1}{-\kappa - \epsilon + \sqrt{\epsilon^2 + \frac{8}{\pi}}})$.\\

Therefore :
\begin{equation*}
\sum_{i=1}^{\infty} \sup_{m \geq n} \log(b_{i,m}) \leq \xi(\kappa, qn) c_{\epsilon}
\end{equation*}

And since $\phi$ is increasing with respect to $p$, we have : 
\begin{equation*}
\limsup_{p \rightarrow \infty} \log(\prod_{i=1}^p \mathbb{P}(\chi_{p/q}^2 \geq \frac{p/q \lambda_{1,p}}{\lambda_{i,p}})) \leq  c_{\epsilon} \lim_{n \rightarrow \infty} \xi(\kappa, qn)
\end{equation*}
Note that $\lim_{n \rightarrow \infty} \phi(\kappa, qn)$ can be finite or infinite. In both cases, and since this is true for any $\epsilon>0$, we have that :
\begin{equation*}
\limsup_{p \rightarrow \infty} \log(\prod_{i=1}^p \mathbb{P}(\chi_{p/q}^2 \geq \frac{p/q \lambda_{1,p}}{\lambda_{i,p}})) \leq  c \lim_{n \rightarrow \infty} \xi(\kappa, qn)
\end{equation*}
where $c = \log(\sqrt{\frac{2}{\pi}} \frac{1}{-\kappa + \sqrt{ \frac{8}{\pi}}}) = \log(\frac{1}{-\kappa \sqrt{\frac{\pi}{2}}  + 2 })$.\\

Since $\kappa < \sqrt{\frac{2}{\pi}}$, then $c<0$, and by taking $c_{\kappa} = - c$, we can conclude.

\section{Interpretation of the results}
In this section, we show in a qualitative way how the conditions of theorem II are reflected on the spectrum when it converges to a limiting spectrum (convergence in terms of empirical density) without outliers (the limiting density has no Dirac mass).\\

Note that it is equivalent to consider the spectrum or the density when the dimension is infinite (we can deduce the spectrum from the density using quantiles, and the density from the spectrum using a histogram).\\
Now Suppose that the spectrum (empirical density function) converges to a deterministic spectrum given by some function $s : x \rightarrow s(x)$ for $x \in [0,1]$ where s is decreasing ($s(0)$ is the biggest eigenvalue and $s(1)$ is the smallest). Let us investigate what the condition $\lim_{p \rightarrow \infty} \phi(p) = \infty$ of theorem II means.\\ 

For this purpose, we define the sequence of spectrums $((\lambda_{i,m})_{1 \leq i \leq m})_{m \geq1}$ by :\\
$\lambda_{i,m} = g(\frac{i}{m})$ where $g : x \rightarrow g(x)$ is a continuous strictly decreasing function defined on $[0,1]$ and differentiable on $]0,1[$.\\

Let p be an positive integer. Recall the definition of the set $J_p$ in theorem II :
\begin{equation}
J_p = \{i : for \hspace{0.1cm} all \hspace{0.1cm} m \geq p, \hspace{0.1cm}  |\frac{\lambda_{1,m}}{\lambda_{i,m}} - 1| < \frac{1}{\sqrt{m}} \}
\end{equation}

We have that :
\begin{align*}
i \in J_p & \Longleftrightarrow \forall m \geq p, \hspace{0.1cm} g(\frac{i}{m}) \geq \frac{g(\frac{1}{m})}{1 + \frac{1}{\sqrt{m}}} \\
& \Longleftrightarrow \forall m \geq p, \hspace{0.1cm} \frac{i}{m} \leq g^{-1}(\frac{g(\frac{1}{m})}{1 + \frac{1}{\sqrt{m}}})
\end{align*}

Since g is a diffeomorphism and g is bounded, we have that :
\begin{align*}
g^{-1}(\frac{g(\frac{1}{m})}{1 + \frac{1}{\sqrt{m}}}) &= g^{-1}(g(\frac{1}{m})(1 - \frac{1}{\sqrt{m}}) +O(\frac{1}{m}))\\
&= g^{-1}(g(\frac{1}{m})) - \frac{g(\frac{1}{m})}{\sqrt{m}} (g^{-1})^{\prime} (g(\frac{1}{m})) +O(\frac{1}{m})\\
&= \frac{1}{m} - \frac{g(\frac{1}{m})}{\sqrt{m}} \frac{1}{g^{\prime}(\frac{1}{m})} + O(\frac{1}{m})\\
\end{align*}

Using this, we have that : 
\begin{align*}
i \in J_p & \Longleftrightarrow for \hspace{0.1cm} all \hspace{0.1cm} m \geq p, \hspace{0.1cm} i \leq 1 - g(\frac{1}{m}) \frac{\sqrt{m}}{g^{\prime}(\frac{1}{m})} + O(1)
\end{align*}

Since $g$ is continuous, we have $\lim_{m \rightarrow \infty} g(\frac{1}{m}) = g(0) > 0$. So in order to have infinite cardinal number of $J_p$ when p goes to infinity, we need to have $\lim_{m \rightarrow \infty} \frac{\sqrt{m}}{g^{\prime}(\frac{1}{m})} = - \infty$, which implies the following property on g : $\lim_{x \rightarrow 0^+}\sqrt{x} g^{\prime}(x) = 0$. From this, we have proved the following proposition.\\

{\bf Proposition 1 }: Let g be a continuous strictly decreasing function defined on $[0,1]$ and differentiable on $]0,1[$. We define the sequence of spectrums $((\lambda_{i,p})_{1 \leq i \leq p})_{p \geq 1}$ for any $p\geq 1$ and $1\leq i \leq p$ by :
\begin{equation}
\lambda_{i,p} = g(\frac{i}{p})
\end{equation}
Let $((l_{i,p})_{1 \leq i \leq p})_{p \geq 1}$ be the corresponding sequence of sample eigenvalues (computed with some given number of samples $n = p/q$ where $q \in ]0,1[ $ fixed).\\

Suppose that $\lim_{x \rightarrow 0^+}\sqrt{x} g^{\prime}(x) = 0$. Then we have that :
\begin{equation*}
\lim_{p \rightarrow \infty} \mathbb{P}(l_{1,p} > \lambda_{1,p}) = 1
\end{equation*}\\

Note that if $\lim_{x \rightarrow 0^+} g^{\prime}(x)$ is finite, then the constraint on g is satisfied. When this limit is infinite (i.e. $\lim_{x \rightarrow 0^+} g^{\prime}(x) = -\infty$ since g is decreasing), the divergence rate should be less that the square root.\\

The following proposition is a general version of the previous result.\\

{\bf Proposition 2 }: Let $((\lambda_{i,p})_{1 \leq i \leq p})_{p \geq 1}$ be a sequence of spectrums and $((l_{i,p})_{1 \leq i \leq p})_{p \geq 1}$ be the corresponding sequence of sample eigenvalues (computed with some given number of samples $n = p/q$ where $q \in ]0,1[ $ fixed). Suppose there exists a function g continuous and strictly decreasing defined on $[0,1]$, differentiable on $]0,1[$ and satisfies $\lim_{x \rightarrow 0^+} \sqrt{x} g(x) = 0$, such that :
\begin{equation*} 
\exists k>0, \forall p \geq k, \forall i \in \{1,...,p\}, \frac{\lambda_{i,p}}{\lambda_{1,p}} \geq \frac{g(\frac{i}{p})}{g(\frac{1}{p})}
\end{equation*}

Then we have that :
\begin{equation*}
\lim_{p \rightarrow \infty} \mathbb{P}(l_{1,p} > \lambda_{1,p}) = 1
\end{equation*}\\

{\it Proof} :\\

We define the sets $J_p$ and $J_p(g)$ by : 
\begin{equation*}
J_p = \{i : for \hspace{0.1cm} all \hspace{0.1cm} m \geq p, \hspace{0.1cm}  |\frac{\lambda_{1,m}}{\lambda_{i,m}} - 1| < \frac{1}{\sqrt{m}} \}
\end{equation*}
and
\begin{equation*}
J_p(g) = \{i : for \hspace{0.1cm} all \hspace{0.1cm} m \geq p, \hspace{0.1cm}  |\frac{g(\frac{1}{m})}{g(\frac{i}{m})} - 1| < \frac{1}{\sqrt{m}} \}
\end{equation*}
It is easy to see that under the assumptions of the proposition we have :
\begin{align*}
J_p(g) \subset J_p
\end{align*}
and using the proof of Proposition 1 we know that $\lim_{p \rightarrow \infty} |J_p(g)| = \infty$, which gives us $\lim_{p \rightarrow \infty} |J_p| = \infty$, we conclude with Theorem II.

\section{A special case : Density with a Dirac mass}
Here we discuss the special case of a density with Dirac mass and see how Theorem II applies in this case.\\

We consider a sequence of spectrums $((\lambda_{i,p})_{1 \leq i \leq p})_{p \geq 1}$ such that, for any $p$, the first $k(p)$ eigenvalues are equal to $\lambda_{1,p}$, where $k(p)$ is a function of p.\\
When the empirical density converges to some limiting density, and $\lim_{p \rightarrow \infty} \frac{k(p)}{p} = \delta \in ]0,1[$ and the rest of the spectrum converges is isolated from the largest eigenvalue, then the limiting density has a Dirac mass on the largest eigenvalue with weight $\delta$.\\
Note that in this case, the conditions of Theorem II are satisfied (since $\phi(p) > k(p) \rightarrow \infty$) and we conclude that with probability 1, we overestimate the largest eigenvalue when the dimension is infinite.\\

Note that the previous argument is true fro any Dirac mass in the limiting density, this proves (using Corollary 1) that the largest sample eigenvalue is bigger than any dirac mass in the limiting distribution when the dimension is infinite.\\

The same argument stands for the smallest eigenvalue using Theorem III. We conclude that we underestimate the smallest eigenvalue when it has a Dirac mass in the limiting density.\\

\section{Conclusion}

\hspace*{0.55cm}In this paper, we proved that for a large class of covariance matrices, we overestimate (underestimate) the largest (smallest) eigenvalue with probability 1 when the dimension is infinite. The conditions on the sets $J_p$ and $H_p$ were derived to be the minimum requirement to have the overestimation (underestimation) using Muirhead's inequality. Thus, the quality of the conditions depends on the quality of the inequalities, and finding more relaxed conditions is directly related to finding better inequalities on the distribution of the extreme eigenvalues of the covariance matrix. 

\section{Acknowledgments}
This work was done during my research internship at Bloomberg LP (New York). I would like to thank Mr Bruno Dupire for giving me the opportunity to work on this project, and I would like to thank all the members of the Quantitative Research team who I worked with.

\newpage


\end{document}